\documentclass[12pt]{amsart}
\usepackage{amssymb,amsmath,amsfonts,latexsym}
\usepackage{bm}
\newcommand{\newsection}[1]{\setcounter{equation}{0} \section{#1}}
\newcommand{\bea}{\begin{eqnarray}}
\newcommand{\eea}{\end{eqnarray}}

\newcommand{\vp}{\varphi}

\newcommand{\clb}{\mathcal{B}}

\newcommand{\cle}{\mathcal{E}}

\newcommand{\clh}{\mathcal{H}}
\newcommand{\clk}{\mathcal{K}}

\newcommand{\clm}{\mathcal{M}}
\newcommand{\cln}{\mathcal{N}}
\newcommand{\clo}{\mathcal{O}}

\newcommand{\cls}{\mathcal{S}}

\newcommand{\z}{\bm{z}}
\newcommand{\w}{\bm{w}}
\newcommand{\K}{\bm{k}}
\newcommand{\AL}{\bm{l}}
\newcommand{\T}{\bm{T}}
\newcommand{\M}{\bm{M}}
\newcommand{\D}{\mathbb{D}}

\newcommand{\raro}{\rightarrow}

\def \qed {\hfill \vrule height6pt width 6pt depth 0pt}
\def\textmatrix#1&#2\\#3&#4\\{\bigl({#1 \atop #3}\ {#2 \atop #4}\bigr)}
\def\dispmatrix#1&#2\\#3&#4\\{\left({#1 \atop #3}\ {#2 \atop #4}\right)}
\newcommand{\be}{\begin{equation}}
\newcommand{\ee}{\end{equation}}
\newcommand{\ben}{\begin{eqnarray*}}
\newcommand{\een}{\end{eqnarray*}}

\newcommand{\NI}{\noindent}

\newcommand{\bi}{\begin{itemize}}
\newcommand{\ei}{\end{itemize}}

\newtheorem{Theorem}{\sc Theorem}[section]
\newtheorem{Lemma}[Theorem]{\sc Lemma}
\newtheorem{Proposition}[Theorem]{\sc Proposition}
\newtheorem{Corollary}[Theorem]{\sc Corollary}

\theoremstyle{definition}

\newtheorem{Remark}[Theorem]{\sc Remark}

\theoremstyle{plain}

\theoremstyle{definition}

\numberwithin{equation}{section}

\let\phi=\varphi

\begin{document}

\title{Factorizations of Kernels and Reproducing Kernel Hilbert Spaces }

\author[Kumari]{Rani Kumari}
\address{Indian Statistical Institute, Statistics and Mathematics Unit, 8th Mile, Mysore Road, Bangalore, 560059, India}
\email{rani\_vs@isibang.ac.in}

\author[Sarkar]{Jaydeb Sarkar}
\address{Indian Statistical Institute, Statistics and Mathematics Unit, 8th Mile, Mysore Road, Bangalore, 560059, India}
\email{jay@isibang.ac.in, jaydeb@gmail.com}

\author[Sarkar]{Srijan Sarkar}
\address{Indian Statistical Institute, Statistics and Mathematics Unit, 8th Mile, Mysore Road, Bangalore, 560059, India}
\email{srijan\_rs@isibang.ac.in}

\author[Timotin] {Dan Timotin}
\address{Institute of Mathematics of the Romanian Academy, P.O. Box 1-764, Bucharest, 014700, Romania }
\email{Dan.Timotin@imar.ro}



\begin{abstract}
The paper discusses a series of results concerning reproducing
kernel Hilbert spaces, related to the factorization of their
kernels. In particular, it is proved that for a large class of
spaces isometric multipliers are trivial. One also gives for certain
spaces conditions for obtaining a particular type of dilation, as
well as a classification of Brehmer type submodules.
\end{abstract}

\maketitle

\newsection{Introduction}

Reproducing kernel Hilbert space theory is an interdisciplinary
subject that arises from the interaction between function theory,
system theory and operator theory. The main aim of this paper is to
investigate the structure of factors of a kernel function and to
relate it with reproducing kernel Hilbert spaces and operators acting
on them.

The precise definition of reproducing kernels is given in Section~\ref{sec:prel}; they may be either scalar or operator valued (the latter type being less familiar to operator theorists).
If $k_1$ is a scalar-valued kernel and $K_2$ is a
$\clb(\cle)$-valued kernel on $\Lambda$, then  $K = K_1 K_2$, where
$K_1 = k_1 I_{\cle}$, is also a
$\clb(\cle)$-valued kernel on $\Lambda$.
We intend to study in the sequel factorizations of reproducing kernels of the
above type and relate function and operator theoretic results on
$\clh_K$ with those of on $\clh_{k_1}$ and $\clh_{K_2}$.

The paper is organized as follows. In Section~\ref{sec:prel} and~\ref{se:kermod} we recall basic
facts concerning reproducing kernel Hilbert spaces, multipliers, and modules over the polynomials. Section~\ref{sec:TPK} is devoted to a presentation of tensor products of reproducing kernel spaces, which are intrinsically related to products of kernels. This part is generally known, but we did not find a suitable reference that would gather all the results we needed.

New results start with Section \ref{sec:IM}, in which we prove that for a large class
of reproducing kernel Hilbert spaces $\clh_K$ isometric multipliers
are trivial. This in particular implies that the reproducing kernel
Hilbert spaces with proper isometrically isomorphic shift invariant
subspaces are rare.

In Section \ref{sec:FKD}, we prove that a reproducing kernel Hilbert
module $\clh_K$ (see the definition in Section \ref{sec:FKD})
defined over a domain $\Omega$ in $\mathbb{C}^n$ dilates to
$\clh_{k_1} \otimes \cle$, for some Hilbert space $\cle$, if and
only if $K = k_1 L$ for some $\clb(\cle)$-valued kernel $L$ on
$\Omega$. Finally, in Section \ref{sec:SRKHS} we obtain a complete
classification of Brehmer type submodules of a large class of
reproducing kernel Hilbert modules and in particular, we prove that
the Brehmer submodules and doubly commuting submodules of the Hardy
module $H^2(\D^n) \otimes \cle$ are the same.

\newsection{Preliminaries}\label{sec:prel}

In this section we briefly recall some basic facts concerning
kernels and reproducing kernel Hilbert spaces. As a general
reference for reproducing kernel Hilbert spaces, see \cite{AM-b} and
\cite{Ar}. For vector-valued reproducing kernel Hilbert spaces,
see~\cite[Chapter 10]{VP}.

Let $\Lambda$ be a set and $\cle$ be a Hilbert space. An
operator-valued function $K : \Lambda \times \Lambda \raro
\clb(\cle)$ is called a \textit{kernel} (cf. \cite{AM-b}, \cite{VP})
and is denoted by $K (\lambda, \mu) \succ 0$, if
\begin{equation}\label{eq:basic formula for definition of kernels}
\sum_{p,q = 1}^m \langle K(x_p, x_q) \eta_q, \eta_p \rangle_{\cle}
\geq 0,
\end{equation}
for all $\{x_j\}_{j=1}^m \subseteq \Lambda$ and $\{\eta_j\}_{j=1}^m
\subseteq \cle$ and $m \in \mathbb{N}$. In this case there exists a
Hilbert space $\clh_K$ of $\cle$-valued functions on $\Lambda$ such
that $\{K(\cdot, \lambda) \eta : \lambda \in \Lambda, \eta \in
\cle\}$ is a total set in $\clh_K$ and
\begin{equation}\label{eq:def reproducing kernel}
\langle f(\lambda), \eta \rangle_{\cle} = \langle f, K(\cdot,
\lambda) \eta \rangle_{\clh_K} \quad \quad (\eta \in \cle, \lambda
\in \Lambda).
\end{equation}
In particular, we have
\begin{equation}\label{eq:norm of reprod kernel}
\|K(\cdot,
\lambda) \eta\|^2_{\clh_K}= \langle K(\lambda, \lambda)\eta, \eta\rangle_\cle=\| K(\lambda, \lambda)^{1/2}\eta\|_\cle.
\end{equation}

\begin{Remark}\label{re:factorization}
    If $\Phi:\Lambda\to \clb(\cle_*, \cle)$ for some Hilbert spaces $\cle,\cle_*$, then it is easy to see that $K(\lambda, \mu):=\Phi(\lambda)\Phi(\mu)^*$ is a kernel with values in $\clb(\cle)$. Conversely, if $K: \Lambda \times \Lambda \raro
    \clb(\cle)$ is a kernel, then we may write $K(\lambda, \mu)=\Phi(\mu)\Phi(\lambda)^*$, with $\cle_*=\clh_K$ and $\Phi(\lambda)=K(\cdot, \lambda)$.
\end{Remark}

Let $\cle$ be a Hilbert space and $K_1$ and $K_2$ be two
$\clb(\cle)$-valued kernel on $\Lambda$. We will write this
sometimes as $K (\lambda, \mu) \succ 0$; then $K_1\prec K_2$ will mean that
 $(K_2-K_1) (\lambda, \mu) \succ 0$.

The following  lemma is known, but for lack of an appropriate reference we supply a
proof for completeness.

\begin{Lemma}\label{le:tensor product of inequalities}
If $K_1 (\lambda, \mu) \prec K_2 (\lambda, \mu)$ and   $L_1
(\lambda, \mu) \prec L_2 (\lambda, \mu)$, then
\[
K_1(\lambda, \mu)\otimes L_1(\lambda, \mu)\prec K_2(\lambda, \mu)\otimes L_2(\lambda, \mu).
\]
\end{Lemma}

\NI \textsf{Proof.} Using~\eqref{eq:basic formula for definition of
kernels}, we have to prove that for nonnegative matrices $A_1, A_2,
B_1, B_2$, if $A_1\le A_2$ and $B_1\le B_2$, then $A_1\otimes B_1\le
A_2\otimes B_2$. One can suppose that $B_1, B_2$ are invertible
(otherwise one adds a small multiple of the identity and pass to the
limit). Therefore
\[
B_1^{-1/2}A_1 B_1^{-1/2}\le I,\quad B_2^{-1/2}A_2 B_2^{-1/2}\le I,
\]
whence (since the tensor product of two contractions is a
contraction)
\[
(B_1^{-1/2}\otimes B_2^{-1/2})(A_1\otimes A_2)( B_1^{-1/2}\otimes
B_2^{-1/2})\le I\otimes I
\]
(the identities acting on the corresponding spaces). It remains to
multiply on the right and on the left with $B_1^{1/2}\otimes
B_2^{1/2}$. \qed

The proof of the following simple lemma is left to the reader.

\begin{Lemma}\label{le:transport of kernel}
Let $K$ be a $\clb(\cle)$-valued kernel on $\Lambda$ and $\clh_\clk$
the corresponding reproducing kernel Hilbert space. Suppose
$\rho:\Lambda'\to\Lambda$ is a bijection. Then
$\clh':=\{f\circ\rho:f\in\clh\}$ endowed with the scalar product
\[
\langle f\circ\rho, g\circ\rho\rangle_{\clh'} :=\langle f,g\rangle_\clh,
\]
is a reproducing kernel Hilbert space of functions on $\Lambda'$,
with the $\clb(\cle)$-valued kernel \[
K'(\lambda',\mu')=K(\rho(\lambda'), \rho(\mu')).
\]\
Moreover, the map $f\mapsto f\circ\rho$ is unitary from $\clh$ to
$\clh'$.
\end{Lemma}

Let $\cle_1$ and $\cle_2$ be two Hilbert spaces and $K_j : \Lambda
\times \Lambda \raro \clb(\cle_j)$, $j=1, 2$, be two kernels. A
function $\varphi : \Lambda \raro \clb(\cle_1, \cle_2)$ is said to
be a \textit{multiplier} if
\[
\varphi f \in \clh_{K_2} \quad {\rm for \; every} \; f \in
\clh_{K_1}.
\]
We will denote by $\clm(\clh_{K_1}, \clh_{K_2})$ the space of all
multipliers from $\clh_{K_1}$ into $\clh_{K_2}$. When $K_1 = K_2$,
we will simply denote it by $\clm(\clh_{K_1})$. From the closed
graph theorem it follows that each multiplier $\varphi \in
\clm(\clh_{K_1}, \clh_{K_2})$ induces a bounded multiplication
operator $M_{\varphi}$ from $\clh_{K_1}$ to $\clh_{K_2}$, where
\[
(M_{\varphi} f)(\lambda) = (\varphi f)(\lambda) = \varphi(\lambda)
f(\lambda) \quad\quad (f \in \clh_{K_1}, \lambda \in \Lambda).
\]
Moreover, for each $\varphi \in \clm(\clh_{K_1}, \clh_{K_2})$,
$\lambda \in \Lambda$ and $\eta \in \cle_2$ we have
\begin{equation}\label{eq:kernel-ev}
M_{\varphi}^* (K_2(\cdot, \lambda) \eta) = K_1(\cdot, \lambda)
\varphi(\lambda)^* \eta.
\end{equation}

We shall call a multiplier $\varphi \in \clm(\clh_{K_1},
\clh_{K_2})$ \textit{partially isometric} or \textit{isometric} if
the induced multiplication operator $M_{\vp}$ has the corresponding
property.

A criterion for multipliers is given in~\cite[Theorem 10.22]{VP}:
$\phi:\Lambda\to\clb(\cle_1, \cle_2)$ is a multiplier if and only if
there exists $c>0$ such that
\begin{equation}\label{eq:criterion for multipliers}
\phi(\lambda)K_1(\lambda, \mu)\phi(\mu)^*\prec c^2 K_2(\lambda, \mu),
\end{equation}
and the smallest such $c$ is precisely the norm of $M_\phi$.

An important particular case are the {\it quasiscalar kernels}. These
are  $\clb(\cle)$-valued kernels of the form
\[
K(\lambda, \mu) = k(\lambda, \mu) I_{\cle} \quad \quad (\lambda,
\mu \in \Lambda),
\]
where $k$ is a scalar-valued kernel on $\Lambda$ and $\cle$ is a
Hilbert space. It follows then from~\eqref{eq:norm of reprod kernel} that
\begin{equation}\label{eq:norm quasi-scalar}
\|K(\cdot, \lambda)\eta\|_{\clh_K} = k(\lambda, \lambda)\|\eta\|_\cle.
\end{equation}
We also note that as Hilbert spaces, one has
\[
\clh_K=\clh_k\otimes \cle.
\]
Therefore, for a fixed orthonormal basis $\{e_j\}$ in $\cle$, the
general form of $F\in \clh_K$ is given by
\[
F=\sum_j f_j\otimes e_j,
\]
with $f_j\in\clh_k$ and $\sum_j\|f_j\|_{\clh_k}^2<\infty$.

Now let $k$ be a scalar kernel and $\lambda \in \Lambda$. By virtue
of (\ref{eq:kernel-ev}), it follows that the functions in $\clh_k$
vanishing at $\lambda$ are given by
\[
\clh_k \ominus \{k(\cdot, \lambda)\} = \{f \in \clh_k : f(\lambda) =
0\}.
\]
For quasiscalar kernels, we have the following:

\begin{Lemma}\label{le:about orthogonals}
Let $k$ be a scalar kernel, $\cle$ a Hilbert  space, $\{e_j\}$ an
orthonormal basis in $\cle$, and $K=kI_{\cle}$ the corresponding
quasiscalar kernel. If $\lambda\in\Lambda$, then
\[
\clh_K\ominus \{k(\cdot, \lambda) x: x\in\cle  \}= \{F=\sum_j
f_j\otimes e_j:f_j\in \clh_k, \  f_j(\lambda)=0, \
\sum_j\|f_j\|_{\clh_k}^2<\infty\}.
\]
\end{Lemma}

\NI\textsf{Proof.} Let us denote by $X$ the space in the right hand
side of the equality. If $F\in X$, then it is immediate that $F$ is
orthogonal to any function $k(\cdot, \lambda) x$.

\NI Conversely, suppose $g=\sum_j g_j \otimes e_j$ is orthogonal to
$X$, that is,
\[
0 = \langle g, F\rangle =\sum_j \langle g_j, f_j\rangle,
\]
for all $F = \sum_j f_j \otimes e_j \in X$. In particular, each
$g_j$ is orthogonal to the space $\{f\in\clh_k: f(\lambda)=0\}$, and
is thus a scalar multiple of $k(\cdot, \lambda) $. Therefore
$g=k(\cdot, \lambda)x$ for some $x\in\cle$. \qed

\section{Kernels and modules}\label{se:kermod}

We now consider a bounded domain $\Omega$ in $\mathbb{C}^n$ and a a
$\clb(\cle)$-valued kernel $K$ on $\Omega$. In what follows, $\z$
will denote the element $(z_1, \ldots, z_n) \in \mathbb{C}^n$.

Let $K(\z, \w)$ be holomorphic in $\{z_1, \ldots, z_n\}$ and
anti-holomorphic in $\{w_1 \ldots, w_n\}$ and $\clh_K$ be the
corresponding reproducing kernel Hilbert space. Then $\clh_K$ is a
set of $\cle$-valued holomorphic functions on $\Omega$ and
\[
\{K(\cdot, \w) \eta : \w \in \Omega, \eta \in \cle\},
\]
is a total set in $\clh_K$, that is,
\[
\clh_K = \overline{\mbox{span}} \{K(\cdot, \w) \eta : \w \in \Omega,
\eta \in \cle\} \subseteq \clo(\Omega, \cle).
\]
In what follows, we always assume that $K(\cdot, \lambda) \neq 0$
for all $\lambda \in \Lambda$.

We say that $\clh_K$ is a \textit{reproducing kernel Hilbert module}
if
\[
z_j \clh_K \subseteq \clh_K \quad \quad (j = 1, \ldots, n).
\]
In this case the multiplication operator tuple $(M_{z_1}, \ldots,
M_{z_n})$, defined by
\[
(M_{z_j} f)(\w) = w_j f(\w) \quad \quad (\w \in \Omega, f \in
\clh_K),
\]
induce a $\mathbb{C}[\z]$-module action on $\clh_K$ as follows:
\[
p \cdot h = p(M_{z_1}, \ldots, M_{z_n}) h \quad \quad (p \in
\mathbb{C}[z_1, \ldots, z_n], h \in \clh_K).
\]
A closed subspace $\cls$ of $\clh_K$ is said to be a
\textit{submodule} if $\cls$ is $M_{z_j}$-invariant, $j = 1, \ldots,
n$. Here the $\mathbb{C}[\z]$-module action on $\cls$ is induced by
the restriction of the multiplication operator tuple $\M_z|_{\cls} =
(M_{z_1}|_{\cls}, \ldots, M_{z_n}|_{\cls})$.

Note also that a submodule of a reproducing kernel Hilbert module is
also a reproducing kernel Hilbert module.

If $\clh_K$ is a reproducing kernel Hilbert module over $\mathbb{C}[\z]$, and the constant functions $\eta\in\cle$ belong to $\clh_K$, then of course  $\mathbb{C}[\z]\cle\subset\clh_K$. The following lemma is often used in concrete cases.

\begin{Lemma}\label{le:intertwining}
Suppose $\clh_{K_i} \subseteq \clo(\Omega, \cle_i)$, $i = 1, 2$ are
reproducing kernel Hilbert modules over $\mathbb{C}[\z]$, and $T:\clh_{K_1}\to \clh_{K_2}$ satisfies
\[
T M_{z_j} = M_{z_j} T \quad \quad (j = 1,
\ldots, n).
\]
If $\mathbb{C}[\z]\cle\subset\clh_K$ and is dense therein, then $T$ is a multiplier.
\end{Lemma}

\NI\textsf{Proof.}
Define $\Phi:\Omega\to\clb(\cle_1, \cle_2)$ by $\Phi(z)\eta=T(\eta)$, where $T(\eta)$ is the action of $T$ on the constant function $z\mapsto\eta\in\cle_1$.
The intertwining assumption in the statement implies that $T(P(z)\eta)=M_\Phi P(z)\eta$ for any polynomial $P$ and $\eta\in\cle$. If $\mathbb{C}[\z]\cle$ is dense in $\clh_K$, it follows that $T=M_{\Phi}$.
\qed

Let $\clh_{K_i} \subseteq \clo(\Omega, \cle_i)$, $i = 1, 2$, be
reproducing kernel Hilbert modules over $\mathbb{C}[\z]$.
We say that they are unitarily equivalent if there exists a unitary $U:\clh_{K_1}\to \clh_{K_2}$ that satisfies
\[
U M_{z_j} = M_{z_j} U \quad \quad (j = 1,
\ldots, n).
\]

\begin{Corollary}\label{co:unitary intertwining}
    With the assumptions of Lemma~\ref{le:intertwining}, $\clh_{K_1}$ and $\clh_{K_2}$ are unitarily equivalent if and only if there exists a unitary multiplier $M_\Phi$ such that $U=M_\Phi$.
\end{Corollary}

\section{Tensor products of kernels}\label{sec:TPK}

Our purpose in this section is to explore the relationship between
kernels and functions defined on a set $\Lambda$ and others defined
on the diagonal of $\Lambda\times\Lambda$.

Let $\cle_i$ are Hilbert spaces and $K_i$ are $\clb(\cle_i)$-valued
kernels on $\Lambda$, $i = 1, 2$. Then the Hilbert tensor product
$\clh_{\clk_1}\otimes \clh_{\clk_2}$ is  a reproducing kernel
Hilbert space on $\Lambda\times\Lambda$, with the
$\clb(\cle_1\otimes\cle_2)$-valued kernel
\[
(K_1\otimes K_2) ((\lambda_1,\lambda_2), (\mu_1,\mu_2)) =
K_1(\lambda_1,\mu_1)\otimes K_2(\lambda_2,\mu_2).
\]
More precisely, the map defined on simple tensors by  $f\otimes
g\mapsto f(\lambda_1)g(\lambda_2)$ extends to a unitary operator
from  $\clh_{\clk_1}\otimes \clh_{\clk_2}$ onto $\clh_{K_1\otimes
K_2}$, which allows the identification of these two spaces.

For clarity, it is useful to make apparent the argument of
functions, typically $\lambda\in\Lambda$ and $(\lambda_1,
\lambda_2)\in \Lambda\times\Lambda$. So, for instance, we will write
$K(\lambda, \mu)$ rather than $K(\cdot, \mu)$ in order to denote the
function $\lambda\mapsto K(\lambda, \mu)$.

Now let $\Delta = \{(\lambda, \lambda): \lambda \in \Lambda\}$ be
the diagonal of $\Lambda \times \Lambda$ and let $\cln$ be the set
of all functions in $ \clh_{K_1}\otimes \clh_{K_2}$ vanishing on
$\Delta$, that is,
\[
\cln=\{g\in \clh_{k_1}\otimes \clh_{k_2}: g(\lambda, \lambda)=0, \
\lambda \in \Lambda\}.
\]
Define also $\delta:\Lambda\to\Delta$ to be the bijection
\[
\delta(\lambda)=(\lambda, \lambda) \quad \quad \quad (\lambda \in
\Lambda).
\]
The scalar case of the next
lemma appears in~\cite{Ar}; we include the proof of the vector case
for completion.

\begin{Lemma}\label{le:cln orthogonal}
With the above notations,
\[
(K_1\ast K_2)(\lambda, \mu):=  K_1(\lambda, \mu) \otimes
K_2(\lambda, \mu).
\]
is a $\clb(\cle_1\otimes\cle_2)$-valued reproducing kernel for the
Hilbert space of functions on $\Lambda$ defined by $\{f\circ\delta :
f\in\cln^\perp \}$, endowed with the scalar product
\[
\langle f\circ\delta, g\circ\delta\rangle_\clh:= \langle
f,g\rangle_{\cln^\perp}.
\]
The map $f\mapsto  f\circ\delta$ is unitary from $\cln^\perp$ to
$\clh_{K_1\ast K_2}$.
\end{Lemma}
\NI\textsf{Proof.} Note first that $\cln^\perp$ is spanned by the
set
\[S:=\{ K_1(\lambda_1, \mu)x_1\otimes K_2(\lambda_2, \mu)x_2 :
\mu\in\Lambda, x_1\in\cle_1, x_2\in\cle_2 \}.
\]
Indeed, for any $F\in\cln$ we have
\begin{equation}\label{eq:about cln}
\langle F, K_1(\lambda_1, \mu)x_1\otimes K_2(\lambda_2, \mu)x_2
\rangle = \langle F(\mu, \mu), x_1\otimes x_2\rangle=0,
\end{equation}
whence $S\subset \cln^\perp$. On the other hand, if $F\in S^\perp$,
then~\eqref{eq:about cln} is true for all $\mu\in\Lambda$ and
$x_1\in\cle_1$, $x_2\in\cle_2$. By linearity we may deduce that
$\langle F(\mu, \mu), \xi\rangle=0$ for all $\xi\in
\cle_1\otimes\cle_2$, whence $F\in \cln$.

It follows then easily that the restrictions of the functions in
$\cln^\perp$ to $\Delta$ form a reproducing kernel Hilbert space,
with kernel given by $K_1(\lambda, \mu) \otimes K_2(\lambda, \mu)$.
The proof is finished by applying Lemma~\ref{le:transport of
kernel}, with $\rho=\delta$. \qed

The proof of the above lemma yields the following useful result:

\begin{Corollary}\label{co:the coisometry pi}
The formula $(\pi F)(\lambda):=F(\lambda, \lambda)$ defines a
coisometry from $\clh_{K_1}\otimes\clh_{K_2}$ to $\clh_{K_1\ast
K_2}$, with  $\ker \pi=\cln$. Also,
\[
\pi^*(  (K_1\ast K_2)(\lambda, \mu)(x_1\otimes x_2)= K_1(\lambda_1,
\mu)x_1 \otimes K_2(\lambda_2, \mu)x_2.
\]
\end{Corollary}
\NI\textsf{Proof.} We observe that, for any $x_j, y_j \in \cle_j$,
$j= 1, 2$ and $\mu, \nu \in \Lambda$,
\[
\begin{split}
\langle K_1(\lambda_1, \mu) x_1 \otimes K_2(\lambda_2, \mu) x_2, &
K_1(\lambda_1, \nu) y_1 \otimes K_2(\lambda_2, \nu) y_2 \rangle =
\langle K_1(\nu, \mu) x_1, y_1 \rangle \langle K_2(\nu, \mu) x_2,
y_2 \rangle
\\
& = \langle K_1(\lambda, \mu) x_1 \otimes K_2(\lambda, \mu) x_2,
K_1(\lambda, \nu) y_1 \otimes K_2(\lambda, \nu) y_2 \rangle.
\end{split}
\]
Then $X : \clh_{K_1\ast K_2} \raro \clh_{K_1} \otimes \clh_{K_2}$
defined by
\[
\pi^*(  (K_1\ast K_2)(\lambda, \mu)(x_1\otimes x_2)= K_1(\lambda_1,
\mu)x_1 \otimes K_2(\lambda_2, \mu)x_2,
\]
for $x_j \in \cle_j$, $j = 1,2$, and $\mu \in \Lambda$, is an
isometry. By the proof of the previous lemma we have $\ker X^* =
(\mbox{ran~} L)^\perp = \cln$ and the result now follows by defining
$\pi = X^*$. \qed

Suppose now that $F_1 : \Lambda \raro \clb(\cle_1)$  is a multiplier
on $\clh_{K_1}$ and $F_2: \Lambda \raro \clb(\cle_2)$ is a
multiplier on $\clh_{K_2}$. Then the  $F_1\otimes
F_2:\Lambda\otimes\Lambda\raro \clb(\cle_1\otimes\cle_2)$ is a
multiplier on $\clh_{K_1\otimes K_2}$, and $M_{F_1\otimes F_2}=
M_{F_1}\otimes M_{F_2}$. The space $\cln$ is invariant to
multipliers on $\clh_{K_1\otimes K_2}$, and therefore $\cln^\perp$
is invariant to adjoints of multipliers.

\begin{Lemma}\label{le:products of multipliers}
If $F_1$ is a multiplier on $\clh_{K_1}$ and $F_2$ is a multiplier
on $\clh_{K_2}$, then the function $F_1\ast F_2:\Lambda\to
\clb(\cle_1\otimes\cle_2) $, defined by $(F_1\ast F_2)(\lambda)=
F_1(\lambda)\otimes F_2(\lambda)$, is a multiplier on $\clh_{K_1\ast
K_2}$. Moreover
\begin{equation}\label{eq:products of multipliers}
M^{K_1\ast K_2}_{F_1\ast F_2}=\pi (M^{K_1}_{F_1}\otimes
M^{K_2}_{F_2}) \pi^*.
\end{equation}
\end{Lemma}
\NI\textsf{Proof.} The assumption implies that~\eqref{eq:criterion
for multipliers} is satisfied for the two multipliers, so
\[
F_1(\lambda)K_1(\lambda, \mu)F_1(\mu)^*\prec c_1^2 K_1(\lambda,
\mu), \quad F_2(\lambda)K_2(\lambda, \mu)F_2(\mu)^*\prec c_2^2
K_2(\lambda, \mu).
\]
By Lemma~\ref{le:tensor product of inequalities}, we have
\[
(F_1(\lambda)\otimes F_2(\lambda))\big( K_1(\lambda, \mu) \otimes
K_2(\lambda, \mu) \big)   (F_1(\mu)\otimes F_2(\mu))^* \prec
c_1^2c_2^2 \big(K_1(\lambda, \mu) \otimes K_2(\lambda, \mu)\big) ,
\]
which means precisely that
\[
(F_1\ast F_2)(\lambda)  (K_1\ast K_2)(\lambda, \mu) (F_1\ast
F_2)(\mu) \prec c_1^2c_2^2 (K_1\ast K_2)(\lambda, \mu).
\]
Again using~\eqref{eq:criterion for multipliers} it follows that
$F_1\ast F_2$ is a multiplier on $\clh_{K_1\ast K_2}$ (of norm at
most $c_1c_2$).

\NI To obtain formula~\eqref{eq:products of multipliers}, we will
check its adjoint on the reproducing kernels $(K_1\ast K_2)(\lambda,
\mu) (x_1\otimes x_2) $, where $\mu\in\Lambda, x_1\in\cle_1,
x_2\in\cle_2$ are fixed, while $\lambda\in\Lambda$ is the variable.
According to ~\eqref{eq:kernel-ev}, we have
\[
(M^{K_1\ast K_2}_{F_1\ast F_2})^* (K_1\ast K_2)(\lambda, \mu) (x_1\otimes x_2) K_1(\lambda, \mu)= K_1(\lambda, \mu)F_1(\mu)^* x_1 \otimes  K_2(\lambda, \mu)F_2(\mu)^* x_1.
\]
On the other hand, by Corollary~\ref{co:the coisometry pi}
\[
\pi^*   (K_1\ast K_2)(\lambda, \mu)(x_1\otimes x_2)= K_1(\lambda_1,
\mu)x_1 \otimes K_2(\lambda_2, \mu)x_2.
\]
Then, applying again~\eqref{eq:kernel-ev},
\[
\big((M^{K_1}_{F_1})^*\otimes (M^{K_2}_{F_2})^*\big)
    \pi^*(  (K_1\ast K_2)(\lambda, \mu)(x_1\otimes x_2)= \big(K_1(\lambda_1, \mu)F_1(\mu)^*x_1\big) \otimes \big(K_2(\lambda_2. \mu)F_2(\mu)^*x_2\big),
\]
Therefore
\[
\pi\big((M^{K_1}_{F_1})^*\otimes (M^{K_2}_{F_2})^*\big)
\pi^*(  (K_1\ast K_2)(\lambda, \mu)(x_1\otimes x_2)=\big( K_1(\lambda, \mu)F_1(\mu)^*x_1\big) \otimes\big( K_2(\lambda, \mu)F_2(\mu)^*x_2\big),
\]
and~\eqref{eq:products of multipliers} is thus proved.
\qed

If one of the kernels is scalar-valued, say  $\dim\cle_2=1$, the
kernel $K_1\ast k_2$ becomes simply the product $k_2K_1$. Then
Lemma~\ref{le:products of multipliers} says that $f_2F_1$ is a
multiplier on $\clh_{k_2K_1}$.

\section{Isometric Multipliers}\label{sec:IM}

In this section, we study the isometric multipliers of reproducing
kernel Hilbert spaces.

Let $k$ be a scalar-valued kernel on a set $\Lambda$ and let
$\clh_k$ be the corresponding reproducing kernel Hilbert space. For
each $\lambda$ and $\mu$ in $\Lambda$, define a relation $\sim_k$ as
follows: $\lambda\sim_k\mu$ if there exist $m \in \mathbb{N}$ and
$\{\lambda_1,\dots,\lambda_m\} \subseteq \Lambda$ such that
\[
\lambda_1=\lambda, \lambda_m=\mu, \mbox{~and~~} k(\lambda_j,
\lambda_{j+1})\not=0 \; \mbox{for~} 1\le j \le m-1.
\]
Then $\sim_k$ is an equivalence relation on $\Lambda$. In
particular, if $\lambda, \mu$ are in two different equivalence
classes, then $k(\lambda, \mu)=0$.

For each multiplier $\varphi \in \clm(\clh_k)$, we will denote,
below, the corresponding multiplication operator $M_\varphi$ by
$M^k_\varphi$.

Suppose $k_1, k_2$ are two scalar-valued reproducing kernels on
$\Lambda$, $K_1=k_1 I_{\cle_1}$, $K_2=k_2 I_{\cle_2}$. If $\phi$ is
a multiplier on $\clh_{K_1}$, it follows from Lemma~\ref{le:products
of multipliers}, applied for $F_1=\phi$ and $F_2=I$, that $\phi*I$
is also a multiplier of $\clh_{K_1*K_2}$, and
\begin{equation}\label{eq:the general norm inequality for multipliers}
\|M^{K_1*K_2}_{\phi*I}\|\le \|M^{K_1}_\phi\|.
\end{equation}

\begin{Theorem}\label{th:general result}
Let  $k_1, k_2$ be two scalar-valued reproducing kernels on
$\Lambda$, and $K_1=k_1 I_{\cle_1}$, $K_2=k_2 I_{\cle_2}$. Denote $k
= k_1 k_2$, $K=K_1*K_2$, and suppose the following conditions are
satisfied:
\begin{enumerate}
\item  the map $M^{K_1}_\phi\mapsto M^{K}_{\phi*I}$  from $\clm(\clh_{K_1})$ to $\clm(\clh_{K})$ preserves the norm;
\item $\clh_{k_1}\cap \clh_{k_2}$ is dense in
$\clh_{k_1}$.
\end{enumerate}
If $M^{K}_{\phi*I}$ is an isometric multiplier  in $\clm(\clh_K)$, then it is a
constant isometry on each of the equivalence classes of~$\sim_{k_1}$.

In particular, if $k_1(\lambda, \mu)\not=0$ for any $\lambda, \mu$,
then $\sim_{k_1}$ has a single equivalence class, and the conclusion
becomes that $\vp$ is a  constant isometry.
\end{Theorem}

\NI\textsf{Proof.}
We use the notation of the previous section; so
$\Delta = \{(\lambda, \lambda): \lambda \in \Lambda\}$, $
\cln=\{F\in \clh_{K_1}\otimes \clh_{K_2}: F(\lambda, \lambda)=0, \
\lambda \in \Lambda\}
$, and $(\pi F)(\lambda)=F(\lambda, \lambda)$ defines a coisometry from $\clh_{K_1}\otimes
\clh_{K_2}$ to $\clh_K$ with $\ker \pi = \cln$.
Then  $M^{K_1}_\phi$ a contraction by assumption (1), and we have by Lemma~\ref{le:products of multipliers}
\begin{equation}\label{eq:commutation of multipliers}
M^K_{\vp\otimes I}= \pi (M^{K_1}_\vp\otimes I_{\clh_{K_2}})\pi^*,
\end{equation}
whence
\[
 M^K_{\vp\otimes I}\pi=\pi  (M^{K_1}_\vp\otimes I_{\clh_{K_2}})\pi^*\pi=\pi  (M^{K_1}_\vp\otimes I_{\clh_{K_2}})P_{\cln^\perp}.
\]
Now for $F\in \cln^\perp$ and using the fact that $M^K_{\vp\otimes I}$ is an
isometry, we have
\[
\|\pi(M^{K_1}_\vp\otimes I_{\clh_{K_2}})F\| =\|M^K_\vp \pi F\|=\|\pi
F\|=\|F\|,
\]
where the last equality follows from the fact that $\pi$ is an
isometry on $(\ker \pi)^\perp$. Hence, since $M^{K_1}_\vp\otimes
I_{k_2}$ is a contraction, we have
\[
\|F\| \geq \|(M^{K_1}_\vp\otimes I_{\clh_{K_2}}) F\| \geq
\|\pi(M^{K_1}_\vp\otimes I_{\clh_{K_2}})F\| = \|F\|,
\]
and hence
\[
\|\pi (M^{K_1}_\vp\otimes I_{\clh_{K_2}}) F\| =
\|(M^{K_1}_\vp\otimes I_{\clh_{K_2}})F\|.
\]
Consequently, $(M^{K_1}_\vp\otimes I_{\clh_{K_2}}) F \in (\ker
\pi)^\perp = \cln^\perp$, that is,
\[
(M^{K_1}_\vp\otimes I_{\clh_{K_2}}) \cln^\perp \subseteq \cln^\perp.
\]
In particular, since $k_1(\lambda_1, \mu)x_1 \otimes k_2(\lambda_2,
\mu)x_2\in\cln^\perp$ for $\mu\in\Lambda$, $x_1\in\cle_1$, $x_2\in\cle_2$ (here $\lambda_1, \lambda_2$ are the argument variables), we have
\[
M^{K_1}_\vp k_1(\lambda_1, \mu)x_1 \otimes k_2(\lambda_2,
\mu)x_2 \in \cln^\perp \quad \quad (\mu \in
\Lambda,\ x_1\in\cle_1, \ x_2\in \cle_2).
\]

Now, if $f,g\in \clh_{k_1}\cap \clh_{k_2}$, $y_1\in\cle_1$, $y_2\in\cle_2$, then $f(\lambda_1)y_1\otimes g(\lambda_2)y_2-g(\lambda_1)y_1
\otimes f(\lambda_2)y_2\in \cln$, and therefore
\[
\begin{split}
0 &= \langle (M^{K_1}_\vp k_1(\lambda_1, \mu)x_1\otimes k_2(\lambda_2,
\mu))x_2, f(\lambda_1)y_1\otimes g(\lambda_2)y_2-g(\lambda_1)y_1
\otimes f(\lambda_2)y_2 \rangle_{\clh_{k_1}
\otimes \clh_{k_2}} \\
&= \langle \vp(\lambda_1) k_1(\lambda_1, \mu)x_1, f(\lambda_1)y_1 \rangle_{\clh_{K_1}} \langle
k_2(\lambda_2, \mu)x_2, g(\lambda_2)y_2 \rangle_{\clh_{K_2}}\\
&\qquad\qquad - \langle \vp(\lambda_1)
k_1(\lambda_1, \mu), g(\lambda_1)y_1 \rangle_{\clh_{K_1}} \langle k_2(\lambda_2, \mu)x_2, f(\lambda_2)y_2 \rangle_{\clh_{K_2}}\\
&  = \langle \vp(\lambda_1) k_1(\lambda_1, \mu)x_1, f(\lambda_1)y_1
\rangle_{\clh_{K_1}} \overline{g(\mu)}\langle x_2, y_2\rangle -\\
&\qquad\qquad
 - \langle \vp(\lambda_1)
k_1(\lambda_1, \mu), g(\lambda_1)y_1 \rangle_{\clh_{K_1}}
\overline{f(\mu)}\langle x_2, y_2\rangle.
\end{split}
\]
Applying assumption (2), the above formula is valid by continuity
for any $f,g\in \clh_{k_1}$.

Fix $\mu\in \Lambda$. Take $f\perp k_1(\lambda, \mu)$ (so
$f(\mu)=0$) and $g=k_1(\cdot, \mu)$ (so $g(\mu)\not=0$); also, assume $\langle x_2, y_2\rangle\not=0$.
It follows from the preceding equation that
\[
\langle \vp(\lambda_1) k_1(\lambda_1, \mu)x_1, f(\lambda_1)y_1 \rangle_{\clh_{K_1}} = 0
\]
for all $x_1, y_1\in\cle_1$. Therefore the function $\vp(\lambda_1)
k_1(\lambda_1, \mu)x_1=M^{k_1}_\vp k_1(\lambda_1, \mu)x_1$ is
orthogonal to the space spanned by the functions $f(\lambda_1)y_1\in
\clh_{K_1}$ with $f\in\clh_{k_1}$, $f(\mu)=0$, and $y_1\in\cle_1$.
If we identify $\clh_{K_1}$ with $\clh_{k_1}\otimes \cle_1$, this
space becomes the space spanned by $f\otimes y_1$, $f(\mu)=0$. We
may then apply Lemma~\ref{le:about orthogonals} to conclude that
\[
\vp(\lambda_1) k_1(\lambda_1, \mu)x_1=k_1(\lambda_1, \mu)x_1'
\]
for some $x_1'\in\cle_1$. But we have, for all $y\in\cle_1$,
\[
\begin{split}
\langle k_1(\lambda_1, \mu)x_1', k_1(\lambda_1, \mu)y\rangle&=
\langle M^{k_1}_\vp
k_1(\lambda_1, \mu)x_1, k_1(\lambda_1, \mu)y\rangle\\
&=\langle
k_1(\lambda_1, \mu)x_1,(M^{k_1}_\vp)^* k_1(\lambda_1, \mu)y\rangle\\
&=\langle
k_1(\lambda_1, \mu)x_1,\vp(\mu)^* k_1(\lambda_1, \mu)y\rangle\\
&=\langle
\vp(\mu)k_1(\lambda_1, \mu)x_1, k_1(\lambda_1, \mu)y\rangle.
\end{split}
\]
Therefore
\[
k_1(\lambda_1, \mu)x_1'=\vp(\lambda_1) k_1(\lambda_1, \mu)x_1
=\vp(\mu) k_1(\lambda_1, \mu)x_1
\]
for all $x_1\in\cle_1$. In this relation $\lambda$ is still the
argument of the functions in the two sides of the equality, but we
may deduce from here the pointwise equality
\[
\vp(\lambda_1) k_1(\lambda_1, \mu)=\vp(\mu) k_1(\lambda_1, \mu).
\]
for all $\lambda_1, \mu\in\Lambda$. So, if $k_1(\lambda_1, \mu)
\not=0$, then $\vp(\lambda_1)=\vp(\mu)$. From the definition of
$\sim_{k_1}$ it follows that on each of its equivalence classes the
multiplier $\phi$ on $\clh_{K_1}$ is a constant operator. Regarding
again $\clh_{K_1}$ as $\clh_{k_1}\otimes \cle_1$, it follows that
$M^{K_1}_\phi=I_{\clh_{k_1}}\otimes \Phi$ for some
$\Phi\in\clb(\cle_1)$. Therefore, in order for $M^{K_1}_\phi$ to be
an isometry, $\Phi$ must be an isometry; this finishes the proof of
the theorem.
\qed

\begin{Corollary}\label{co:of the main theorem}
Let  $k_1, k_2$ be two scalar-valued reproducing kernels on
$\Lambda$, and $K_1=k_1 I_{\cle_1}$. Denote  $K=K_1k_2$, and suppose
the following conditions are satisfied:
\begin{enumerate}
\item  the map $M^{K_1}_\phi\mapsto M^{K}_{\phi}$  from $\clm(\clh_{K_1})$ to $\clm(\clh_{K})$ is surjective and preserves the norm;
\item $\clh_{k_1}\cap \clh_{k_2}$ is dense in $\clh_{k_1}$.
\end{enumerate}
Then any isometric multiplier  in $\clm(\clh_K)$ is a constant
isometry on each of the equivalence classes of~$\sim_{k_1}$.
\end{Corollary}

There is an important case in which condition (1) in the above
corollary is satisfied, which we will present as a separate
statement.

\begin{Corollary}\label{co:case of good multiplier norm}
Let $\Lambda = \Omega$ be a domain in $\mathbb{C}^n$ and $k_1, k_2$
are analytic in the first variable, $K_1=k_1 I_{\cle_1}$,
$K=K_1k_2$. Suppose the following conditions are satisfied:
\begin{enumerate}
\item
$\clm(\clh_{K_1})$ coincides with the uniformly bounded
$\clb(\cle_1)$-valued analytic functions and for any
$\phi\in\clm(\clh_{K_1})$ we have
\begin{equation}\label{eq:equality for norm of multipliers}
\|M_\phi\|_{\clh_{K_1}}=    \sup_\lambda
\|\phi(\lambda)\|;\end{equation}
\item $\clh_{k_1}\cap \clh_{k_2}$ is dense in $\clh_{k_1}$.
\end{enumerate}
Then any isometric multiplier  in $\clm(\clh_K)$ is a
constant isometry on each of the equivalence classes of~$\sim_{k_1}$.
\end{Corollary}
\NI\textsf{Proof.} Let $\phi\in \clm(\clh_K)$. Then by
(\ref{eq:kernel-ev}), we have
\[
\sup_\lambda \|\phi(\lambda)\|\le \|M_\phi^K\|_{\clb(\clh_{K})}.
\]
Hence condition (1) imply that $M^{K_1}_\phi\in \clm(\clh_{K_1})$.
Applying also~\eqref{eq:the general norm inequality for
multipliers}, it follows that $\|M^{K_1}_\phi\|=\|M^{K}_\phi\|$. We
may then apply Corollary~\ref{co:of the main theorem} to conclude
the proof. \qed

\begin{Remark}\label{re:kernel} Under the same assumptions and notations as in
Corollary~\ref{co:case of good multiplier norm}, suppose also that
polynomials are in $\clh_{k_1}$ as well as in $\clh_{k_2}$. Then a
sufficient condition for (2) is that they are dense in $\clh_{k_1}$.
\end{Remark}

Using now Corollary~\ref{co:unitary intertwining}, it is easy to derive the following result.

\begin{Theorem}\label{th:unit equiv}
    Let $\cle$ be a Hilbert space, $\Omega$ be a domain in
    $\mathbb{C}^n$ and $\clh_{k_1}$, $\clh_{k_2} \subseteq \clo(\Omega)$
    are reproducing kernel Hilbert spaces. Let $K_1 = k_1 I_{\cle}$ and
    $K = k_1 k_2 I_{\cle}$. Suppose the following conditions are
    satisfied:

    \begin{enumerate}
        \item $\clh_{k_1}$ is a reproducing kernel Hilbert module over $\mathbb{C}[\z]$.

        \item $\mathbb{C}[\z] \subseteq \clh_{k_1} \cap \clh_{k_2}$ and
        $\mathbb{C}[\z]$ is dense in $\clh_{k_1}$.

        \item  $\clm(\clh_{K_1}) =
        H^\infty_{\clb(\cle)}(\Omega)$ and for each $\varphi
        \in\clm(\clh_{K_1})$ we have
        \[
        \|M_\phi\|_{\clh_{K_1}}=    \sup_{\lambda \in \Omega}
        \|\phi(\lambda)\|.
        \]
        \item $z_1\sim_{k_1} z_2$ for any $z_1, z_2\in\Omega$. In particular, this is true if $k_1(z_1, z_2)\not =0$ for any  $z_1, z_2\in\Omega$.
    \end{enumerate}
    Let $\varphi \in\clm(\clh_{K_1})$ be a multiplier and $\cls$ be a
    submodule of $\clh_K$. Then

    (i) $M_{\varphi}$ is an isometric multiplier if and only if there
    exists an isometry $V \in \clb(\cle)$ such that $M_{\varphi} =
    I_{\clh_{k_1 k_2}} \otimes V$.

    (ii) $\cls \subseteq \clh_K$ is unitarily equivalent to $\clh_K$ if
    and only if there exists a closed subspace $\tilde{\cle}$ of $\cle$
    such that $\cls = \clh_{k_1 k_2} \otimes \tilde{\cle}$.
\end{Theorem}

Let $\alpha >n$ and $g(\z, \w) = (1 - \sum_{i=1}^n z_i
\bar{w}_i)^{-\alpha}$, $\z , \w \in \mathbb{B}^n$. Then $\clh_g$,
also denoted by $L^2_{a, \alpha}(\mathbb{B}^n)$, is a weighted
Bergman module over $\mathbb{B}^n$. It is well known that the
multiplier space of $L^2_{a, \alpha}(\mathbb{B}^n)$ is
$H^\infty(\mathbb{B}^n)$.

The following corollary is now immediate:

\begin{Corollary}
    Let $\cle$ be a Hilbert space $\varphi \in
    H^\infty_{\clb(\cle)}(\mathbb{B}^n)$ and $\cls$ be a submodule of
    $L^2_{a, \alpha}(\mathbb{B}^n) \otimes \cle$. Then

    (i) $M_{\varphi}$ is an isometry if and only if  $M_{\varphi} =
    I_{L^2_{a, \alpha}(\mathbb{B}^n)} \otimes V$ for some sometry $V \in
    \clb(\cle)$.

    (ii) $\cls$ be a unitarily equivalent to $L^2_{a,
        \alpha}(\mathbb{B}^n) \otimes \cle$ if and only if $\cls = L^2_{a,
        \alpha}(\mathbb{B}^n) \otimes \tilde{\cle}$ for some closed subspace
    $\tilde{\cle}$ of $\cle$.
\end{Corollary}

Part (ii) of the above theorem is related to the rigidity of
submodules of weighted Bergman modules (see \cite{DS-1, GHX, P, R}).
Part (i) is a generalization of Proposition 4.2 in \cite{O1}.

\section{Factorizations of Kernels and dilations}\label{sec:FKD}

A scalar-valued kernel $g$ on $\Omega$ is said to be \textit{good
kernel} if $\clh_g$ is a reproducing kernel Hilbert module and
\[
\mathop{\cap}_{j=1}^n \ker (M_{z_j}^* - \bar{w}_j I_{\clh_g}) =
\mathbb{C} g(\cdot, \w) \quad \quad (\w \in \Omega),
\]
and there exists a $\w_0 \in \Omega$ such that
\[
g(\cdot, \w_0) \equiv 1.
\]
We say that $\clh_g \subseteq \clo(\Omega, \mathbb{C})$ is a
\textit{good reproducing kernel Hilbert module}.

We notice that if $g$ is a scalar valued kernel on a set $\Lambda$
and the function $g(\cdot, \lambda_0)$ is non-vanishing for some
$\lambda_0 \in \Lambda$ then one can assume, after renormalizing,
that $g(\cdot, \lambda_0)\equiv 1$.

Let $\clh_g$ be a good reproducing kernel Hilbert module over
$\Omega$ and $\clh_K \subseteq \clo(\Omega, \cle)$ be a reproducing
kernel Hilbert module over $\mathbb{C}[\z]$. We say that $\M_z =
(M_{z_1}, \ldots, M_{z_n})$ on $\clh_K$ \textit{dilates} to
$(M_{z_1} \otimes I_{\cle}, \ldots, M_{z_n} \otimes I_{\cle})$ on
$\clh_g \otimes \cle$, or $\clh_K$ dilates to $\clh_g \otimes \cle$,
for some Hilbert space $\cle$, if there exists an isometry $\Pi :
\clh \raro \clh_g \otimes \cle$ such that
\[
(M_{z_i}^* \otimes I_{\cle}) \Pi =  \Pi M_{z_i}^* \quad \quad (i =
1, \ldots, n).
\]

Our main result in this section is the following  theorem which
relates dilation of a reproducing kernel Hilbert module to a good
reproducing Hilbert module with factorizations of kernels.

\begin{Theorem}\label{thm-factor}
Let $\cle$ and $\cle_*$ be two Hilbert spaces and $\clh_g$ be a good
reproducing kernel Hilbert module on $\Omega$ and $\clh_K \subseteq
\clo(\Omega, \cle)$ be a reproducing kernel Hilbert module over
$\mathbb{C}[\z]$. Then the following conditions are equivalent:
\begin{enumerate}
\item $\clh_K$ dilates to $\clh_g \otimes \cle_*$.

\item There exists a holomorphic function $\Phi : \Omega \raro
\clb(\cle_*, \cle)$ such that
\[
K(\z, \w) = g(\z, \w) \Phi(z) \Phi(\w)^* \quad \quad (\z, \w \in
\Omega).
\]
\end{enumerate}
\end{Theorem}

\NI \textsf{Proof.}
%
Assume (3) holds. Then for each
$\z, \w \in \Omega$ and $\eta, \zeta \in \cle_*$, we have
\[
\begin{split}
\langle K(\cdot, \w) \eta, K(\cdot, \z) \zeta \rangle_{\clh_K} & =
\langle K(\z, \w) \eta, \zeta \rangle_{\cle_*}
\\
& = \langle g(\z, \w) \Phi(\z) \Phi(\w)^* \eta, \zeta
\rangle_{\cle_*}
\\
& = g(\z, \w) \langle \Phi(\z) \Phi(\w)^* \eta, \zeta
\rangle_{\cle_*}
\\
& = \langle g(\cdot, \w), g(\cdot, \z) \rangle_{\clh_g} \langle
\Phi(\w)^* \eta, \Phi(\z)^* \zeta \rangle_{\cle_*}
\\
& = \langle g(\cdot, \w) \otimes \Phi(\w)^* \eta, g(\cdot, \z)
\otimes \Phi(\z)^* \zeta \rangle_{\clh_g \otimes \cle_*}.
\end{split}
\]
This allows us to define an isometry $\Pi : \clh_K \raro \clh_g
\otimes \cle_*$ by
\[
\Pi(K(\cdot, \w) \eta) = g(\cdot, \w) \otimes \Phi(\w)^* \eta \quad
\quad (\w \in \Omega, \eta \in \cle_*).
\]
Using this, on one hand, we have
\[
\begin{split}
(\Pi M_{z_j}^*)(K(\cdot, \w) \eta) & = \Pi (\bar{w}_j K(\cdot, \w)
\eta)
\\
& = \bar{w}_j \Pi (K(\cdot, \w) \eta)
\\
& = \bar{w}_j (g(\cdot, \w) \otimes \Phi(\w)^* \eta),
\end{split}
\]
and on the other hand, by (\ref{eq:kernel-ev}) , we have
\[
\begin{split}
(M_{z_j} \otimes I_{\cle_*})^* \Pi (K(\cdot, \w) \eta) & = (M_{z_j}
\otimes I_{\cle_*})^* (g(\cdot, \w) \otimes \Phi(\w)^* \eta)
\\
& = \bar{w}_j (g(\cdot, \w) \otimes \Phi(\w)^* \eta),
\end{split}
\]
where $\eta \in \cle$ and $\w \in \Omega$. Therefore
\begin{equation}\label{eq:M-inter}
(M_{z_j} \otimes I_{\cle_*})^* \Pi = \Pi M_j^* \quad \quad (j = 1,
\ldots, n),
\end{equation}
and hence $\clh_K$ dilates to $\clh_g \otimes \cle_*$. This proves
(1).

Assume now (1) hold. Then there exists an isometry $\Pi : \clh_K
\raro \clh_g \otimes \cle_*$ such that (\ref{eq:M-inter}) hold. Then
for $\w \in \Omega$ and $\eta \in \cle$ and $j = 1, \ldots, n$, we
have
\[
\begin{split}
(M_{z_j} \otimes I_{\cle_*})^* (\Pi K(\cdot, \w) \eta) & = ((M_{z_j}
\otimes I_{\cle_*})^* \Pi) (K(\cdot, \w) \eta)
\\
& = \Pi M_{z_j}^* (K(\cdot, \w) \eta)
\\
& = \bar{w}_j (\Pi K(\cdot, \w) \eta).
\end{split}
\]
In particular,
\[
\Pi (K(\cdot, \w) \eta) \in \mathop{\cap}_{j=1}^n \ker \Big((M_{z_j}
\otimes I_{\cle_*})^* - \bar{w}_j I_{\clh_g \otimes \cle_*}\Big) =
g(\cdot, \w) \otimes \cle_*.
\]
Then for each $\w \in \Omega$ there exists a linear map $\Phi(\w) :
\cle_* \raro \cle$ such that
\[
\Pi(K(\cdot, \w) \eta) = g(\cdot, \w) \otimes \Phi(\w)^* \eta \quad
\quad (\eta \in \cle).
\]
Observe that if $\w \in \Omega$ and $\eta \in \cle$ we have
\[
\begin{split}
\|\Phi(\w)^* \eta\|_{\cle_*} & = \frac{1}{\| g(\cdot, \w)
\|_{\clh_g}} \|\Pi(K(\cdot, \w) \eta)\|_{\clh_g \otimes \cle_*}
\\
& \leq \frac{1}{\| g(\cdot, \w) \|_{\clh_g}} \|(K(\cdot, \w)
\eta)\|_{\clh_K}
\\
& \leq \frac{1}{\| g(\cdot, \w) \|_{\clh_g}} \|K(\w,
\w)^{\frac{1}{2}}\|_{\clb(\cle)} \|\eta\|_{\cle},
\end{split}
\]
where the last inequality follows from the fact that
\[
\begin{split}
\|(K(\cdot, \w) \eta)\|^2_{\clh_K} & = \langle K(\cdot, \w) \eta,
K(\cdot, \w) \eta \rangle_{\clh_K} \quad \quad (\mbox{by}
(\ref{eq:def reproducing kernel})) \\ & = \langle K(\w, \w) \eta,
\eta \rangle_{\cle} \\ & = \|K(\w, \w)^{\frac{1}{2}}
\eta\|_{\cle}^2.
\end{split}
\]
Therefore $\Phi(\w)^*$, $\w \in \Omega$, is a bounded linear
operator. For $\eta, \zeta \in \cle$ we now have
\[
\begin{split}
\langle K(\z, \w) \eta, \zeta \rangle_{\cle} & = \langle K(\cdot,
\w) \eta, K(\cdot, \z) \zeta \rangle_{\clh_K}
\\
& = \langle \Pi(K(\cdot, \w) \eta), \Pi (K(\cdot, \z) \zeta)
\rangle_{\clh_g \otimes \cle_*}
\\
& = \langle g(\cdot, \w) \otimes \Phi(\w)^* \eta, g(\cdot, \z)
\otimes \Phi(\z)^* \zeta \rangle_{\clh_g \otimes \cle_*}
\\
& = g(\z, \w) \langle \Phi(\w)^* \eta, \Phi(\z)^* \zeta
\rangle_{\cle_*}
\\
& = \langle g(\z, \w) \Phi(\z) \Phi(\w)^* \eta, \zeta
\rangle_{\cle_*},
\end{split}
\]
and hence
\[
K(\z, \w) = g(\z, \w) \Phi(\z) \Phi(\w)^* \quad \quad (\z, \w \in
\Omega).
 \]
Finally, since
\[
\begin{split}
\langle \Phi(\w) \zeta, \eta \rangle_{\cle} & = \langle \zeta,
\Phi(\w)^* \eta \rangle_{\cle_*}
\\
& = \langle g(\cdot, \w_0) \otimes \zeta, g(\cdot, \w) \otimes
\Phi(\w)^* \eta \rangle_{\clh_g \otimes \cle_*}
\\
& = \langle g(\cdot, \w_0) \otimes \zeta, \Pi (K(\cdot, \w)\eta)
\rangle_{\clh_g \otimes \cle_*}
\\
& = \langle \Pi^* (g(\cdot, \w_0) \otimes \zeta) ,K(\cdot, \w)\eta
\rangle_{\clh_g \otimes \cle_*},
\end{split}
\]
for each $\eta \in \cle$ and $\zeta \in \cle_*$, and since $\w
\mapsto K(\cdot, \w)$ is anti-holomorphic, we conclude that $\w
\mapsto \Phi(\w)$ is holomorphic. This shows that (3) holds and
completes the proof of the theorem. \qed

The next corollary follows by taking into account Remark~\ref{re:factorization}.

\begin{Corollary}\label{co:factorization}
Let $\cle$  be a Hilbert spaces and $\clh_g$ be a good
reproducing kernel Hilbert module on $\Omega$ and $\clh_K \subseteq
\clo(\Omega, \cle)$ be a reproducing kernel Hilbert module over
$\mathbb{C}[\z]$. Then the following conditions are equivalent:
\begin{enumerate}
    \item There exists a Hilbert space $\cle_*$ such that the equivalent conditions in the statement of Theorem~\ref{thm-factor} hold.
    \item There exists a $\clb(\cle)$-valued kernel $L$ on $\Omega$,
    holomorphic in the first and anti-holomorphic in the second
    variable,  such that $K = g L$.
\end{enumerate}
\end{Corollary}

Theorem~\ref{thm-factor} and Corollary~\ref{co:factorization} represent a generalization of the dilation results of quasi-free
Hilbert modules (see Theorems 1 and 2 in \cite{DMS}) to reproducing
kernel Hilbert modules. Let us also note that, moreover, our argument does not rely on
localizations of Hilbert modules.

\section{submodules of reproducing kernel Hilbert
modules}\label{sec:SRKHS}

Let $p(\z,\w) = \sum_{\K, \AL \in \mathbb{N}^n} a_{{\bm{k}} \bm{l}}
\z^{\bm{k}} \bar{\bm{w}}^{\bm{l}}$ be a polynomial in $(z_1, \ldots,
z_n)$ and $(\bar{w}_1, \ldots, \bar{w}_n)$. Here $(z_1, \ldots,
z_n)$ and $(\bar{w}_1, \ldots, \bar{w}_n)$ are commuting variables
but we do not assume commutativity of $z_i$ and $\bar{w}_j$, $1 \leq
i, j \leq n$. Then for a commuting tuple $\T = (T_1, \ldots, T_n)$
on a Hilbert space $\clh$, we define $p(\T, \T^*)$ by
\[
 p(\T, \T^*) = \sum_{\K, \AL \in \mathbb{N}^n} a_{{\bm{k}} \bm{l}} \T^{\bm{k}}
{\T}^{*\bm{l}}.
\]

We will often consider in this section a good kernel $g$ with the property that $g^{-1}$ is a polynomial. We will then write
\[
g^{-1}(\z, \w) = \sum_{\K, \AL \in \mathbb{N}^n} a_{{\bm{k}} \bm{l}}
\z^{\bm{k}} \bar{\bm{w}}^{\bm{l}},
\]
having always in mind that the sum is finite.

The following standard relationship between factorized kernels and
operator positivity of multiplication tuples on reproducing kernel
Hilbert modules is well known (cf. Theorem 4 in \cite{DMS}).

\begin{Proposition}\label{prop-positive}
Let $\clh_K \subseteq \clo(\Omega, \cle)$ be a reproducing kernel
Hilbert module and $g$ be a good kernel on $\Omega$ with $g^{-1}$ a polynomial
Then $g^{-1}(\bm{M}_z, \bm{M}_z) \geq 0$ on $\clh_K$ if and only if
there exists a kernel $L$ on $\Omega$ such that $K = g L$.
\end{Proposition}
\NI \textsf{Proof.} It is enough to prove that $g^{-1}(\M_z, \M_z^*)
\geq 0$ if and only if $g^{-1} K$ is positive definite. Indeed, for
$\{\w_j\}_{j=1}^m \subseteq \Omega$, $\{\eta_j\}_{j=1}^m \subseteq
\cle$ and $m \in \mathbb{N}$, we have
\[
\begin{split}
\sum_{i, j=1}^m \langle (g^{-1} K)(\w_i, \w_j)\eta_j, \eta_i \rangle
& = \sum_{i, j=1}^m \langle g^{-1}(\w_i, \w_j) K(\w_i, \w_j) \eta_j,
\eta_i \rangle
\\
& = \sum_{i, j=1}^m \sum_{\K, \AL \in \mathbb{N}^n} a_{{\bm{k}}
\bm{l}} \w^{\bm{k}} \bar{\bm{w}}^{\bm{l}} \langle  K(\cdot, \w_j)
\eta_j, K(\cdot, \w_i) \eta_i \rangle
\\
& = \sum_{i, j=1}^m \sum_{\K, \AL \in \mathbb{N}^n} a_{{\bm{k}}
\bm{l}} \langle  \M_z^{*\AL}  K(\cdot, \w_j) \eta_j, \M_z^{*\K}
K(\cdot, \w_i) \eta_i \rangle
\\
& = \sum_{i, j=1}^m \langle (\sum_{\K, \AL \in \mathbb{N}^n}
a_{{\bm{k}} \bm{l}} \M_z^{\K} \M_z^{*\AL})  K(\cdot, \w_j) \eta_j,
K(\cdot, \w_i) \eta_i \rangle
\\
& = \sum_{i, j=1}^m \langle g^{-1}(\M_z, \M_z^*)  K(\cdot, \w_j)
\eta_j, K(\cdot, \w_i) \eta_i \rangle.
\end{split}
\]
This completes the proof. \qed

This and Theorem \ref{thm-factor} immediately yields the following
generalization of Theorem 6 in \cite{DMS}.

\begin{Theorem}\label{thm-dil-fac-pos}
In the setting of Proposition \ref{prop-positive} the operator
$g^{-1}(\M_z, \M_z^*) \geq 0$ on $\clh_K$ if and only if there
exists a kernel $L$ on $\Omega$ such that $K = g L$, if and only if
there exists a Hilbert space $\cle_*$ such that $\clh_K$ dilates to
$\clh_g \otimes \cle_*$.
\end{Theorem}

We now turn to the study of submodules of good reproducing kernel
Hilbert modules. To this end, we first need the following simple
lemma.

\begin{Lemma}\label{lem-const}
Let $\clh_g$ on $\Omega$ be a good reproducing kernel Hilbert module
over $\mathbb{C}[\z]$, with $ g^{-1}(\z, \w) = \sum_{\K, \AL
    \in \mathbb{N}^n} a_{{\bm{k}} \bm{l}} \z^{\bm{k}}
\bar{\bm{w}}^{\bm{l}}$ a polynomial. Let $P_{g(\cdot, \w_0)}$ be the orthogonal
projection of $\clh_g$ onto the one dimensional subspace generated
by $g(\cdot, \w_0) \equiv 1$. Then
\[\sum_{\K, \AL \in \mathbb{N}^n} a_{\K \AL} \M_z^{\K} \M_z^{*
\AL} = P_{g(\cdot, \w_0)}.
\]
\end{Lemma}
\NI\textsf{Proof. } For each $\z, \w \in \Omega$ we compute
\[
\begin{split}
\langle \sum_{\K, \AL \in \mathbb{N}^n} a_{\K \AL} \M_z^{\K} \M_z^{*
\AL} g(\cdot, \w), g(\cdot, \z)\rangle & = \sum_{\K, \AL \in
\mathbb{N}^n} a_{\K \AL} \langle \M_z^{\K} \M_z^{* \AL} G(\cdot,
\w), g(\cdot, \z) \rangle
\\
 & = \sum_{\K, \AL \in
\mathbb{N}^n} a_{\K \AL} \langle \M_z^{* \AL} g(\cdot, \w), \M_z^{*
\K} g(\cdot, \z) \rangle
\\
& = (\sum_{i,j=0}^k \z^{\K} \bar{\w}^{\AL}  a_{\K \AL}) \langle
g(\cdot, \w), g(\cdot, \z) \rangle
\\
& = g^{-1} (\z, \w) g(\z, \w) = 1
\\
& = \langle P_{g(\cdot, \w_0)} g(\cdot, \w), g(\cdot, \z)\rangle.
\end{split}
\]
This completes the proof of the lemma. \qed

Let $\clh_g$ be as in the previous lemma and $\cle$ be a Hilbert
space. Let $\cls$ be a submodule of $\clh_g \otimes \cle$, that is,
$\cls$ is a joint $(M_{z_1} \otimes I_{\cle}, \ldots, M_{z_n}
\otimes I_{\cle})$ invariant subspace of $\clh_g \otimes \cle$. Here
$\cls$ is a module over $\mathbb{C}[\z]$ with module multiplication
operators $\bm{R}_z = (R_{z_1}, \ldots, R_{z_n})$, where
\[
R_{z_i} = M_{z_i}|_{\cls} \quad \quad (i = 1, \ldots, n).
\]
We say that $\cls$ is a \textit{Brehmer submodule} if
\[g^{-1}(\bm{R}_z, \bm{R}_z^*) =  \sum_{\K, \AL \in \mathbb{N}^n}
a_{\K \AL} \bm{R}_z^{\K} \bm{R_z}^{* \AL} \geq 0.
\]

In the following we characterize Brehmer submodules in terms of
partial isometric multipliers. The idea of the proof is to invoke
the dilation result, Theorem \ref{thm-dil-fac-pos}, to submodules of
good reproducing kernel Hilbert modules (cf. \cite{JS1}).

\begin{Theorem}\label{thm-breh-sub}
Let $\cle$ be a Hilbert space and $g$ be a good kernel with
$g^{-1}$  a polynomial. Let $\cls$ be a
submodule of $\clh_g \otimes \cle$. Then $\cls$ is a Brehmer
submodule of $\clh_g \otimes \cle$ if and only if there exists a
Hilbert space $\cle_*$ and a partial isometric multiplier $\Theta
\in \clm(\clh_g \otimes \cle_*, \clh_g \otimes \cle)$ such that
\[
\cls = \Theta (\clh_g \otimes \cle_*).
\]
\end{Theorem}
\NI\textsf{Proof.} Let $\cls$ be a Brehmer submodule, that is,
\[
g^{-1}(\bm{R}_z, \bm{R}_z^*) \geq 0.
\]
Then by Theorem \ref{thm-dil-fac-pos}, there exists a Hilbert space
$\cle_*$ such that $\cls$ dilates to $\clh_g \otimes \cle_*$.
Therefore there exists an isometry $\pi : \cls \raro \clh_g \otimes
\cle_*$ such that
\[
\pi R_{z_i}^* = (M_{z_i} \otimes I_{\cle_*})^* \pi \quad\quad (i =
1, \ldots, n).
\]
Let $i : \cls \raro \clh_g \otimes \cle$ be the inclusion map and
$\Pi = i \circ \pi^*$. Then $\Pi : \clh_g \otimes \cle_* \raro
\clh_g \otimes \cle$ is a partial isometry and
\[
\mbox{ran~} \Pi = \cls,
\]
and
\[
\Pi (M_{z_i} \otimes I_{\cle_*}) = (M_{z_i} \otimes I_{\cle}) \Pi
\quad \quad (i = 1, \ldots, n).
\]
This yields that $\Pi = M_{\Theta}$ for some partial isometric
multiplier $\Theta  \clm(\clh_g \otimes \cle_*, \clh_g \otimes
\cle)$ and $\cls = \Theta (\clh_g \otimes \cle_*)$.

Conversely, let $\cls = \Theta (\clh_g \otimes \cle_*)$ for some
partial isometric multiplier $\Theta \in \clm(\clh_g \otimes \cle_*,
\clh_g \otimes \cle)$. Then
\[P_{\cls} = M_{\Theta} M_{\Theta}^*.
\]
Hence
\[
\begin{split}
g^{-1}(\bm{R}_z, \bm{R}_z^*) & =  \sum_{\K, \AL \in \mathbb{N}^n}
a_{\K \AL} \bm{R}_z^{\K} \bm{R_z}^{* \AL}
\\
& =  \sum_{\K, \AL \in \mathbb{N}^n} a_{\K \AL} \M_z^{\K} P_{\cls}
\M_z^{* \AL}
\\
& =  \sum_{\K, \AL \in \mathbb{N}^n} a_{\K \AL} \M_z^{\K} M_{\Theta}
M_{\Theta}^* \M_z^{* \AL}
\\
& =  M_{\Theta}(\sum_{\K, \AL \in \mathbb{N}^n} a_{\K \AL} \M_z^{\K}
\M_z^{* \AL})M_{\Theta}^*
\\
& =  M_{\Theta} P_{g(\cdot, \w_0)} M_{\Theta}^* \quad \quad \quad
(\mbox{by Lemma~} \ref{lem-const})
\\
& \geq 0.
\end{split}
\]
This completes the proof of the theorem.

\qed

Now we consider the important case when $\Omega = \D^n$ and $\clh_g =
H^2(\D^n)$ and $n \geq 2$. A submodule $\cls$ of $H^2(\D^n) \otimes
\cle$ is said to be \textit{doubly commuting} (cf. \cite{SSW}) if
\[
[R_{z_i}^*, R_{z_j}] := R_{z_i}^* R_{z_j} - R_{z_j} R_{z_i}^* = 0,
\]
for all $1 \leq i \neq j \leq n$.

The next theorem is proved in~\cite{SSW}.

  \begin{Theorem}\label{le:doubly commuting}
    A submodule $\cls$ of
    $H^2(\D^n) \otimes \cle$ is doubly commuting if and only if there
    exists a Hilbert space $\cle_*$ and an inner multiplier $\Theta \in
    \clm(H^2(\D^n) \otimes \cle_*, H^2(\D^n) \otimes \cle) =
    H^\infty_{\clb(\cle_*, \cle)}(\D^n)$ such that
    \[
    \cls = \Theta (H^2(\D^n) \otimes \cle_*).
    \]
  \end{Theorem}

In the following, we prove that the class of doubly commuting
submodules and the class of Brehmer submodules of $H^2(\D^n) \otimes
\cle$ are the same.

\begin{Theorem}\label{th:pi=iso}
Let $\cle$ be a Hilbert space. Then  $\cls$ is a Brehmer submodule of
$H^2(\D^n) \otimes \cle$ if and only if $\cls$ is a doubly commuting
submodule.
\end{Theorem}

\NI\textsf{Proof.} If $\cls$ is a doubly commuting submodule, it
follows from Theorem ~\ref{le:doubly commuting} and
Theorem~\ref{thm-breh-sub} that it is a Brehmer submodule.

Conversely, supose $\cls $ is a Brehmer submodule.
By Theorem \ref{thm-breh-sub}, there exists a
Hilbert space $\cle_*$ and a partial isometry $M_{\Theta} :H^2(D^n)
\otimes \cle_* \raro H^2(\D^n) \otimes \cle$, for some multiplier
$\Theta \in H^\infty_{\clb(\cle_*, \cle)}(\D^n)$, such that
\[
\cls = \Theta (H^2(\D^n) \otimes \cle_*).
\]
It is easy to see that the closed subspace $\ker M_{\Theta}$ is a
submodule of $H^2(\D^n) \otimes \cle_*$. We claim that the
orthogonal of $\ker M_{\Theta}$ is also a submodule of $H^2(\D^n)
\otimes \cle_*$. Indeed, if $f \in (\ker M_{\Theta})^\perp$, then
\[
\|f\|=\|M_{z_i} M_{\Theta} f\|=\|M_{\Theta} M_{z_i} f\|\le
\|M_{z_i}f\|=\|f\|,
\]
and hence the inequality becomes an equality. But then
\[
\|M_{\Theta} M_{z_i} f\| = \|M_{z_i} f\|,
\]
yields $z_i f \in (\ker M_{\Theta})^\perp$ for all $i = 1, \ldots,
n$, and hence $(\ker M_{\Theta})^\perp$ is a submodule of $H^2(\D^n)
\otimes \cle_*$, or equivalently that $(\ker M_{\Theta})^\perp$ is a
joint $(M_{z_1}, \ldots, M_{z_n})$-reducing subspace of $H^2(\D^n)
\otimes \cle_*$. Now the reducing subspaces of $H^2(\D^n) \otimes
\cle_*$ are known to be of the form $H^2(\D^n) \otimes
\tilde{\cle}_*$ for some $\tilde{\cle}_* \subseteq \cle_*$. Then
$\cls$ is the image of the isometric multiplier
$M_{\Theta}|_{H^2(\D^n) \otimes \tilde{\cle}_*}$, so $\cls$ is
doubly commuting by the result quoted above. This completes the
proof of the theorem. \qed

\vspace{0.2in}

\NI\textit{Acknowledgement:} The research of the second author was
supported in part by NBHM (National Board of Higher Mathematics,
India) Research Grant NBHM/R.P.64/2014. The research of the fourth
author was supported in part by a grant of the Romanian National
Authority for Scientific Research, CNCS--UEFISCDI, project number
PN-II-ID-PCE-2011-3-0119. The fourth author is also grateful for
hospitality of Indian Statistical Institute, Bangalore, during
November-December 2015.


\begin{thebibliography}{99}

\bibitem{AM-b}
J. Agler and J. McCarthy, {\em Pick Interpolation and Hilbert
Function Spaces}, Grad. Stud. Math., 44. Amer. Math. Soc.,
Providence, RI, 2002.

\bibitem{Ar}
N. Aronszajn, {\em Theory of reproducing kernels}, Transactions of
the American Mathematical Society 68 (1950), 337–-404.

\bibitem{CV}
R. Curto and F.-H. Vasilescu, {\em Standard operator models in the
polydisc}, Indiana Univ. Math. J. 42 (1993), no. 3, 791–-810.

\bibitem{DS-1}
R. G. Douglas and J. Sarkar, {\em On unitarily equivalent
submodules}, Indiana University Mathematics Journal, 57, (2008)
2729--2743.


\bibitem{DMS}
R. G. Douglas, G. Misra and J. Sarkar, {\em Contractive Hilbert
modules and their dilations}, Israel Journal of Math. 187 (2012),
141--165.

\bibitem{GHX}
K. Guo, J. Hu, and X. Xu, {\em Toeplitz algebras, subnormal tuples
and rigidity on reproducing $\mathbb{C}[z_1,\ldots,z_d]$-modules},
J. Funct. Anal. 210 (2004), 214–-247

\bibitem{NF}
B. Sz.-Nagy and C. Foias, {\em Harmonic Analysis of Operators on
Hilbert Space. North-Holland}, Amsterdam-London, 1970.

\bibitem{vN}
J. von Neumann, {\em Eine Spektraltheorie fur allgemeine Operatoren
eines unitdren Raumes}, Math. Nachr. 4 (1951), 258-–281.

\bibitem{O1}
O. Giselsson and A. Olofsson, {\em On some Bergman shift operators},
Complex Anal. Oper. Theory 6 (2012) 829--842.

\bibitem{VP}
V. Paulsen and M. Raghupathi, \emph{An introduction to the theory of
reproducing kernel Hilbert spaces},  Cambridge Studies in Advanced
Mathematics, 2016.

\bibitem{P}
M. Putinar, {\em On the rigidity of Bergman submodules}, Amer. J.
Math. 116 (1994), 1421–-1432.

\bibitem{R}
S. Richter, {\em Unitary equivalence of invariant subspaces of
Bergman and Dirichlet spaces}, Pac. J. Math. 133 (1988), 151–-156

\bibitem{JS1}
J. Sarkar, {\em An invariant subspace theorem and invariant
subspaces of analytic reproducing kernel Hilbert spaces. I,} J.
Operator Theory 73 (2015), 433–-441.


\bibitem{SSW}
J.  Sarkar, A. Sasane, and B. Wick, {\em Doubly commuting submodules
of the Hardy module over polydiscs}, Studia Math. 217 (2013),
179–-192.


\end{thebibliography}
\end{document}